# Adaptive Model Predictive Control of a Batch Solution Polymerization Process using Trajectory Linearization

Masoud Abbaszadeh

*Department of Electrical and Computer Engineering,
University of Alberta, Edmonton, Alberta, Canada, T6G 2V4,*
E-mail: masoud@ualberta.net

**Abstract:** A sequential trajectory linearized adaptive model based predictive controller is designed using the DMC algorithm to control the temperature of a batch MMA polymerization process. Using the mechanistic model of the polymerization, a parametric transfer function is derived to relate the reactor temperature to the power of the heaters. Then, a multiple model predictive control approach is taken in to track a desired temperature trajectory. The coefficients of the multiple transfer functions are calculated along the selected temperature trajectory by sequential linearization and the model is validated experimentally. The controller performance is studied on a small scale batch reactor.

## 1. Introduction

The importance of effective polymer reactor control has been emphasized in recent decades. Kinetic studies are usually complex because of the nonlinearity of the process. Hence, the control of the polymerization reactor has always been a challenging task. Due to its great flexibility, a batch reactor is suitable to produce small amounts of special polymers and copolymers. The batch reactor is always dynamic by its nature. A good dynamic response over the entire process is necessary to reach an effective controller performance. To do so, it is essential to have a suitable dynamic model of the process. Louie *et al* [1] reviewed the gel effect models and their theoretical foundations. These researchers then modeled the solution polymerization of methyl methacrylate (MMA) and validated their model.

Control of MMA polymerization processes has become popular as a benchmark for advaced process control methods, since the dynamics of the methyl methacrylate polymerization process is well studied and several physical models of high fedility are readily avaliable. Methyl methacrylate is normally produced by a free radical, chain addition polymerization. Free radical polymerization consists of three main reactions: initiation, propagation and termination. Free radicals are formed by the decomposition of initiators. Once formed, these radicals propagate by reacting with surrounding monomers to produce long polymer chains; the active site being shifted to the end of the chain when a new monomer is added. Rafizadeh [2] presented a review on the proposed models and suggested an on-line estimation of some parameters, such as heat transfer coefficients. The model consists of the oil bath, electrical heaters, cooling water coil, and reactor. Mendoza-Bustos *et al* [3] derived a first order plus dead time transfer function for polymerization. Then, they designed PID, Smith predictor, and Dahlin controllers for temperature control. Peterson *et al* [4] presented a non-linear predictive strategy for semi batch polymerization of MMA. Penlidis *et al* [5] presented an excellent paper, in which they reviewed a mechanistic model for bulk and solution free radical polymerization for control purposes. Soroush and Kravaris [6] applied a Global Linearizing Control (GLC) method to control the reactor temperature. They compared the result of GLC and PID controllers. Performance of the GLC for tracking an optimum temperature trajectory was found to be suitable. DeSouza jr. *et al* [7] studied an expert neural network as an internal model in control of solution polymerization of vinyl estate. The architecture of their model predicts one step ahead. In their study, they compared their neural network control with a classic PID controller. Clarke-Pringle and MacGregor [8] studied the temperature control of a semi-batch industrial reactor. They suggested a coupled non-linear strategy and extended Kalman filter method. They used energy balance approach for the reactor and jacket to estimate process parameters. Mutha *et al* [9] suggested a non-linear model based control strategy, which includes a new estimator as well as Kalman filter. They conducted experiments in a small reactor for solution polymerization of MMA. Rho *et al* [10] reviewed the batch polymerization modeling and estimated the model parameters based on the experimental data in the literature. For control purposes, they assumed a model to pursue the control studies and estimated the parameters of this model by on line ARMAX model.

Model predictive control (MPC), on the other hand, is a model based advanced control technique that have been proved to be very sussefull in controlling highly complex dynamic systems. It naturally supports design for MIMO and time-delayed systems as well as state/input/output constants. MPC is generally based on online optimization but in the case of unconstrianed linear plants, closed form solutions can be derived analytically. MPC usally requires a high computaional power; however, since chemical processes are typically of slow dynamics, they have been designed and implemented on various chemical plnat with great success. Therefore, MPC seems to be good candidate for controlling MMA polymerization based on physical (first-principle) modeling.

This paper presents a mechanistic model of batch polymerization. Sequential linearization, along a selected temperature trajectory, is conducted. Consequently, using a nonlinear model predictive approach, a controller is designed. A multiple model adaptive MPC controller is desined for the trajectory lineairzed model. Results show the better performance than the performance of adaptive PI controller [11]. Our approach is based on an output feedback control archtecture using (noisy) output measurements. The controller essentially comprise both an implicit state observer and a control law. Since, out approach uses an online trajectory liniearization, the MPC law is a linear controller. Alternative MBC approach based on nonlinear plant model would require an EKF [12] or robust nonlinear obserever such as those proposed in [13-20].

## 2. Polymerization Mechanism

Methyl methacrylate normally is produced by a free radical, chain addition polymerization. Free radical polymerization consists of three main reactions: initiation, propagation and termination. Free radicals are formed by the decomposition of initiators. Once formed, these radicals propagate by reacting with surrounding monomers to produce long polymer chains; the active site being shifted to the end of the chain when a new monomer is added. During the propagation, millions of monomers are added to $P_1^o$ radicals. During termination, due to reactions among free radicals, the concentration of radicals decreases. Termination is by combination or disproportionation reactions. With chain transfer reactions to monomer, initiator, solvent, or even polymer, the active free radicals are converted to dead polymer [1]. Table 1 gives the basic free radical polymerization mechanism.

Table 1. Polymerization mechanism [11]

| | | |
|---|---|---|
| Initiation | $I \xrightarrow{k_d} 2R^o + G^\uparrow$ | $R_d = k_d[I]$ |
| | $R^o + M \xrightarrow{k_i} P_1^o$ | $R_i = k_i[R^o][M]$ |
| | $2R^o \xrightarrow{k_{ti}} I'$ | $R_{ti} = k_{ti}[R^o]^2$ |
| Propagation | $P_n^o + M \xrightarrow{k_p} P_{n+1}^o$ | $R_p = k_p[M][P_n^o]$ |
| Termination | $P_n^o + P_m^o \xrightarrow{k_{tc}} D_{n+m}$ | $R_{tc} = k_{tc}[P_n^o][P_m^o]$ |
| | $P_n^o + P_m^o \xrightarrow{k_{td}} D_n + D_m$ | $R_{td} = k_{td}[P_n^o][P_m^o]$ |
| Transfer | $P_n^o + M \xrightarrow{k_f} P_1^o + D_n$ | $R_f = k_f[M][P_n^o]$ |
| | $P_n^o + M \xrightarrow{k_s} s^o + D_n$ | $R_s = k_s[M][P_n^o]$ |

The free radical polymerization rate decreases due to reduction of monomer and initiator concentration. However, due to viscosity increase beyond a certain conversion there is a sudden increase in the polymerization rate. This effect is called Trommsdorff, gel, or auto-acceleration effect. For bulk polymerization of methyl methacrylate beyond the 20% conversion, reaction rate and molecular weight suddenly increase. In high conversion, because of viscosity increase there is a reduction in termination reaction rate.

## 3. Mathematical Modeling of Polymerization

Table 2 shows the mass and energy balances of reactor. The polymer production is accomplished by a reduction in volume of the mixture. The volumetric reduction factor is given by:

(1) $$\varepsilon = \frac{\rho_p - \rho_m}{\rho_p}$$

The instantaneous volume of mixture is given by:

(2) $$V = \frac{M_0}{\rho_m}(1 - \varepsilon x + \beta)$$

The parameter $\beta$ is defined as:

(3) $$\beta = \frac{f_s}{1 - f_s}$$

Table 2- Mass and energy balances

$$\frac{dx}{dt} = \frac{2 f k_d}{[M]_0 V_0}[I] + (k_p + k_f)(1-x)\lambda_0$$

$$\frac{d[I]}{dt} = -k_d[I] - \frac{[I]}{V}\frac{dV}{dt}$$

$$\frac{d[S]}{dt} = -k_s[S]\lambda_0 - \frac{[S]}{V}\frac{dV}{dt}$$

$$\frac{dV}{dt} = -\frac{M_0}{\rho_m}\left[\varepsilon\frac{dx}{dt} + x\frac{d\varepsilon}{dt}\right]$$

$$mC_p\frac{dT}{dt} = (-\Delta H_p)k_p[M]\lambda_0 V - UA|_r(T - T_j) - UA|_\infty(T - T_\infty)$$

$$m_o C_{p_o}\frac{dT_j}{dt} = \alpha P + \alpha P + UA|_r(T - T_j) - UA|_{o\infty}(T_j - T_\infty)$$

$$\lambda_0 = \sqrt{\frac{2 f k_d [I]}{k_t}}$$

During the free radical polymerization, the cage, glass, and gel effects occur. For the cage effect, the initiator efficiency factor is used. The CCS (Chiu, Carrat, and Soong) model is used in this study to take into consideration the glass and the gel effects. Therefore, propagation rate constant, $k_p$, is changing according to:

(4) $$\frac{1}{k_p} = \frac{1}{k_{p_0}} + \theta_p \frac{\lambda_0}{D}$$

$k_{p_0}$ is changing as Arrhenius function, and $D$ is given by equation:

(5) $$D = \exp\left[\frac{1 - \varphi_p}{A + B(1 - \varphi_p)}\right]$$

Similarly, termination rate constant, $k_t$, is given by:

(6) $$\frac{1}{k_t} = \frac{1}{k_{t_0}} + \theta_t \frac{\lambda_0}{D}$$

$k_{t_0}$ is changing as Arrhenius function. $\theta_p$ and $\theta_t$ are adjustable parameters related to propagation and termination rate constants, respectively. All the other necessary parameters and constants for this model are given in the literature [1]. The equations 7 to 10 are essential for dynamic studies.

(7) $$\frac{dx}{dt} = (k_p + k_f)(1-x)\sqrt{\frac{2fk_d[I]}{k_t}} = f_1(x,[I],T)$$

(8) $$\frac{d[I]}{dt} = -k_d[I] + \frac{\varepsilon}{1-\varepsilon x + \beta}(k_p + k_f)[I](1-x)\sqrt{\frac{2fk_d[I]}{k_t}} = f_2(x,[I],T)$$

(9) $$\frac{dT}{dt} = \frac{(-\Delta H_p)k_p V_0 [M]_0 (1-x)\sqrt{\frac{2fk_d[I]}{k_t}} - UA|_r(T-T_j) - UA|_\infty(T-T_\infty)}{mC_p} = f_3(x,[I],T,T_j)$$

(10) $$\frac{dT_j}{dt} = \frac{2\alpha P + UA|_r(T-T_j) - UA|_{o\infty}(T_j - T_\infty)}{m_o C_{p_o}} = f_4(T,T_j,P)$$

Equations 7 and 8 are mass balances for monomer and initiator, respectively. Long Chain Approximation (LCA) and Quasi Steady State Approximation (QSSA) are used in this study. Equations 9 and 10 show energy balances for the reactant mixture and oil, respectively. In this study, heat transfer coefficients are estimated experimentally [2]. Equations 7 to 10 are highly nonlinear and, using Taylor expansion series, these equations were converted to linearized form. The linearized state space form is given by:

(11) $$\begin{bmatrix} \frac{dX}{dt} \\ \frac{di}{dt} \\ \frac{dT'}{dt} \\ \frac{dT'_j}{dt} \end{bmatrix} = \begin{bmatrix} \frac{\partial f_1}{\partial x}\bigg|_s & \frac{\partial f_1}{\partial [I]}\bigg|_s & \frac{\partial f_1}{\partial T}\bigg|_s & 0 \\ \frac{\partial f_2}{\partial x}\bigg|_s & \frac{\partial f_2}{\partial [I]}\bigg|_s & \frac{\partial f_{21}}{\partial T}\bigg|_s & 0 \\ \frac{\partial f_3}{\partial x}\bigg|_s & \frac{\partial f_3}{\partial [I]}\bigg|_s & -\frac{UA|_r + UA|_\infty}{mC_p} & \frac{UA|_r}{mC_p} \\ 0 & 0 & \frac{UA|_r}{m_o C_{po}} & -\frac{UA|_r + UA|_{\infty o}}{m_o C_{po}} \end{bmatrix} \begin{bmatrix} X \\ i \\ T' \\ T'_j \end{bmatrix} + \begin{bmatrix} 0 \\ 0 \\ 0 \\ \frac{2\alpha}{m_o C_{po}} \end{bmatrix} P'$$

$$T' = \begin{bmatrix} 0 & 0 & 1 & 0 \end{bmatrix} \begin{bmatrix} X \\ i \\ T' \\ T'_j \end{bmatrix}$$

where

(12) $$X = (x - x_s), \quad i = ([I] - [I]_s), \quad T' = (T - T_s), \quad T'_j = (T_j - T_{js}), \quad P' = P - P_s$$

equation 11 and is converted to the transfer function form:

(13) $$\frac{T'(s)}{P(s)} = \frac{n_3 s^2 + n_4 s + n_5}{d_1 s^4 + d_2 s^3 + d_3 s^2 + d_4 s + d_5}$$

## 4. The Experimental Setup

A schematic representation of the experimental batch reactor setup is shown in Figure 1. The reactor is a Buchi type jacketed, cylindrical glass vessel. A multi-paddle agitator mixes the content. A Pentium II 500 MHz computer is connected to the reactor via an ADCPWM-01 analog/digital Input/Output data acquisition card. The data acquisition software was developed in-house. The heating oil was circulated by a gear pump and its flow rate was about $15 \, lit/\min$. The heating/cooling system of the oil consisted of two 1500W electrical heaters and a coolant water coil, which was operated by an On/Off Acco brand solenoid valve. Two Resistance Temperature Detectors (RTDs), were used with accuracy of $\pm 0.2 \, ^oC$. Methyl methacrylate and toluene were used as monomer and solvent, respectively. Benzoyl peroxide (BPO) was used as the initiator. The molecular weight of the produced polymer was measured using an Ubbelohde viscometer.

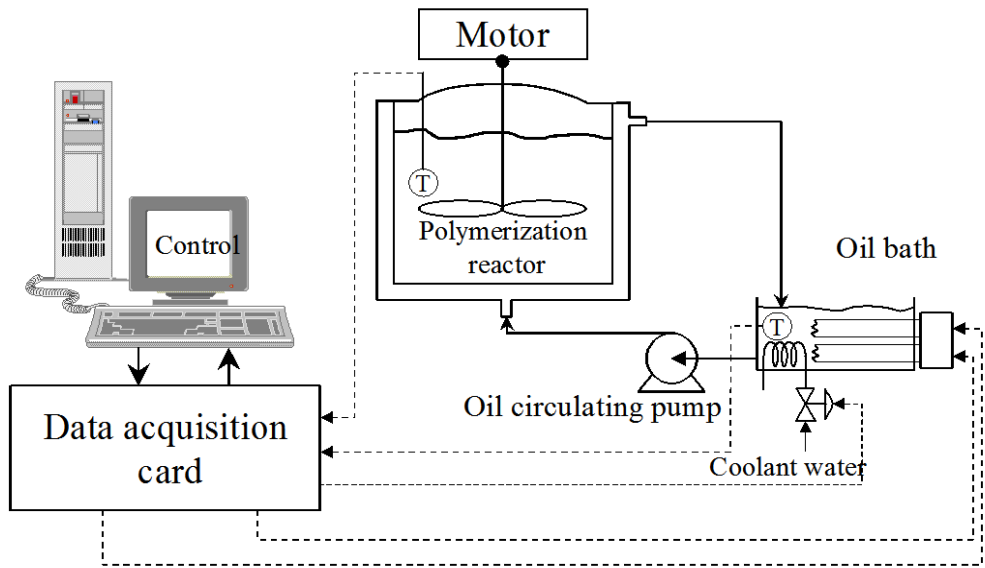

Figure 1. The experimental setup

## 5. An Overview of MPC

Due to its high performance, Model predictive control method has recieved a great deal of attention to control chemical processes, in last few years. This approach is applicable to multivariable systems and canstrained systems. Monuverability in design, noise and disturbance rejection and robustness under model mismatch are the most important ability of this method. Cumbersome computation, lack of systematic rules for controller tuning are some drawback of this method. Model predictive control is based on a process model. Although impulse or step responses have some limitation for nonlinear process, they may be used to develop a model. During the the model predictive control following steps should be conducted:
- Explicit prediction of future output (prediction horizon).
- Calculation of a control sequence based on the minimized cost function (control horizon).
- Receding strategy.

The Dynamic Matrix Control (DMC) is used in this research. Its cost function is:

$$(14) \quad J = \sum_{i=N_1+1}^{N_1+P} \|e(t+i)\|_Q^2 + \sum_{i=1}^{M} \|\Delta u(t+i-1)\|_R^2$$

where P, M and $N_1$ are prediction horizon, control horizon and pure time delay, respectivly. $R_{M \times M}$, $Q_{P \times P}$ are whigthing martices. The prediction horizon must be at least equal to the pure time delay.

$$(15) \quad e(t+i) = y_d(t+i) - y_p(t+i) = y_d(t+i) - y_m(t+i) - d(t+i)$$

where $y_p$ is the process output, $y_m$ is the model output and d is the process and model outputs diffrence, including noise, disturbance and model mismatch. $y_d(t)$ is the desired output based on the refrence input. If $y_{sp}(t)$ is the refrence input, the following filtered form is used as the tracking trajectory:

$$(16) \quad y_d(t) = \alpha y_d(t-1) + (1-\alpha) y_{sp}(t) \quad ; \quad 0 \le \alpha < 1$$

$\alpha$ changes the first order smoothing filter pole place. The smaller $\alpha$ the faster output. It has been showen that system robustness can be decreased by the reduction of $\alpha$ and increment of the manipulated signal [21-22]. Figure 2 shows the block diagram of DMC.

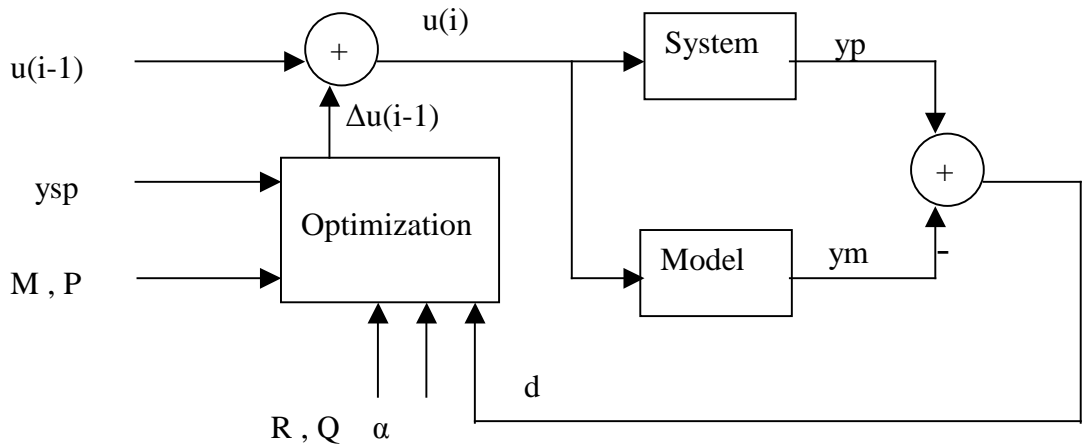

Figure 2. Block diagram of DMC

The cost function in equation 10 can be rearrenged to:

$$(17) \quad J = (Y_m + D - Y_D)^T Q (Y_m + D - Y_D) + \Delta U^T R \Delta U$$

without loss of generality, if $N_1$ is assumed zero, then:

$$Y_m = [y_m(t+1)...y_m(t+P)]^T$$
$$Y_D = [y_d(t+1)...y_d(t+P)]^T$$
$$\Delta U = [\Delta u(t)...\Delta u(t+M-1)]^T \;,\; D = [d(t+1)...d(t+P)]^T$$

For LTI system, without any constraints on output or control signal, optimization has the following closed form:

$$(18) \quad \Delta U_+ = (G_+^T Q G_+ + R)^{-1} G_+^T Q E$$

where:

(19) $$E \overset{\Delta}{=} Y_P - Y_D$$

$$G_+ = \begin{bmatrix} g_\circ & \cdots & \cdots & \cdots \\ g_1 & g_\circ & & \\ \vdots & & \ddots & \\ g_{P-1} & & & \ddots \end{bmatrix}, \quad \Delta U_+ = \Delta U = [\Delta u(t) \cdots \Delta u(t+M-1)]^T$$

$g_i$ s are the step response samples and $G_+$ is a Toeplitz matrix cocsisting the step response samples. The model output has calculated by:

(20) $\quad Y_m = G_+ \Delta U_+ + G_- \Delta U_- + g_N U_N$

$$G_- = \begin{bmatrix} g_1 & g_2 & \cdots & g_{N-1} \\ g_2 & \cdots & g_{N-1} & \circ \\ \vdots & & & \vdots \\ g_p & \cdots & \cdots & \circ \end{bmatrix}, \quad \Delta U_- = [\Delta u(t-1) \cdots \Delta u(t-N+1)]^T$$

where N is the number of system step response samples reaching to steady state or equvalent impulse response steps which lead to zero; and $g_N$ is the system dc gain.

$$g_N = dcgain$$
$$U_N = [u(t-N) \; u(t-N+1) \cdots u(t-N+P-1)]^T$$

## 9. The Modified DMC

If there is any pole close to origin, the step response will be very slow and the required N is very large. Then, a system including integrator never reaches to the steady state (this case exists in the set of linearized models of the MMA reactor) and N lead to infinity. Hence, unstability occurs. This is one of the DMC limitations [23].

Reserchers have suggetsted some methods to overcome this problem, for example formulating DMC in the state space an then using an state observer [24]. Because of model mismatch this method doesn't have proper performance in real time applications. The alternative is:

$$Y_m = G_+ \Delta U_+ + G_- \Delta U_- + g_N U_N$$

(21) $\quad Y_{Past} \overset{\Delta}{=} G_- \Delta U_- + g_N U_N \Rightarrow Y_m = G_+ \Delta U_+ + Y_{Past}$

where $Y_{Past}$ is " the effect of past input to the future system outputs without considering the effect of present and future inputs". Consequently, $Y_{Past}$ can be calculated by setting the future "$\Delta u$"s equal to zero and solving the model P steps ahead.

(22) $\quad \Delta U_+ = \circ \rightarrow Y_m = Y_{Past}$

As seen in equations 14 and 15, $G_+$ and $\Delta U_+$ are independent of N. $G_-$ dimension is determined by N. Therefore, the DMC calculation is independent than N. $Y_D$ is:

(23) $\quad Y_D = [y_d(t+1) \cdots y_d(t+p)]^T$

## 10. Simulation Results

Figure 3 [11], shows the model validation results. The simulation follows the experimental data very well. The DMC algorithm was applied to control a MMA polymerization reactor. The reactor temperature trajectory is known, hence, the refrence input is known for all times so the programmed MPC is used. The DMC controller gain defined as:

(24) $\quad K_{DMC} = (G_+^T Q G_+ + R)^{-1} G_+^T Q$

(25) $\quad \Delta U_+ = k_{DMC} E$

The $G_+$ of present model is used to calculate its $k_{DMC}$. $k_{DMC}$ is changed in the appropriate model switching instant. Therefore, a multiple model stategy is used implicitly. However the valid model is known before.

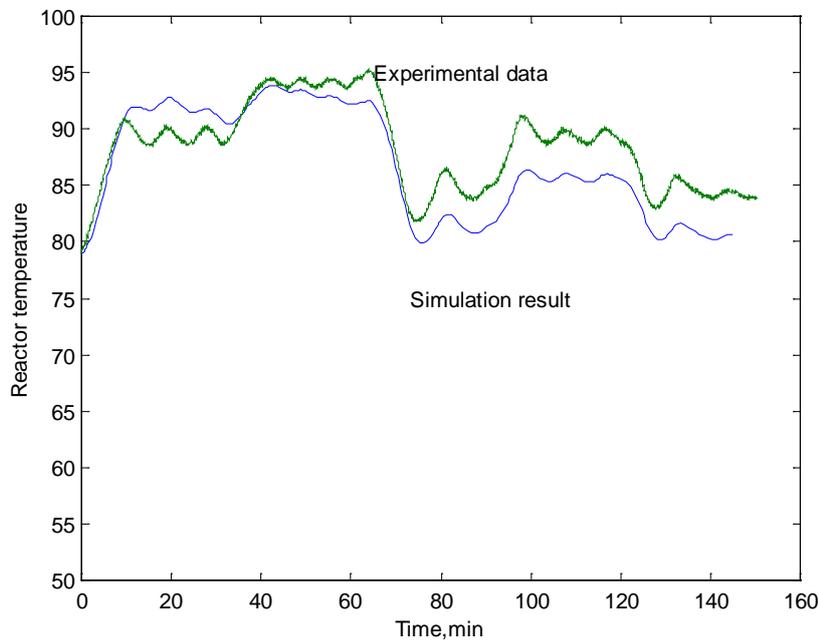

Figure 3. Model validation

Figures 4 and 5 show the controller ability to track the temperature trajectory. The average error is $0.3^oC$. Due to the controller robustness, switching between models causes no unstability in closed loop system. Furthermore, appropriate selection of controller parameters could prevent the unstability. The selected sampling period is T=10s. Other parameters are $P=5$, $M=2$, $\alpha=.05$, $Q=I_{5*5}$, $R=.05*I_{3*3}$. The adaptive multiple model MPC designed here ensures the reactor temperature tracking error to with in $0.3^oC$ while the adaptive PI control in [11] has a *2° C* avarage error and the Generalized Takagi-Sugeno-Kang fuzzy controller proposed in [25] has a *1° C* avarage error; demonestrating the superior performance of the MPC.

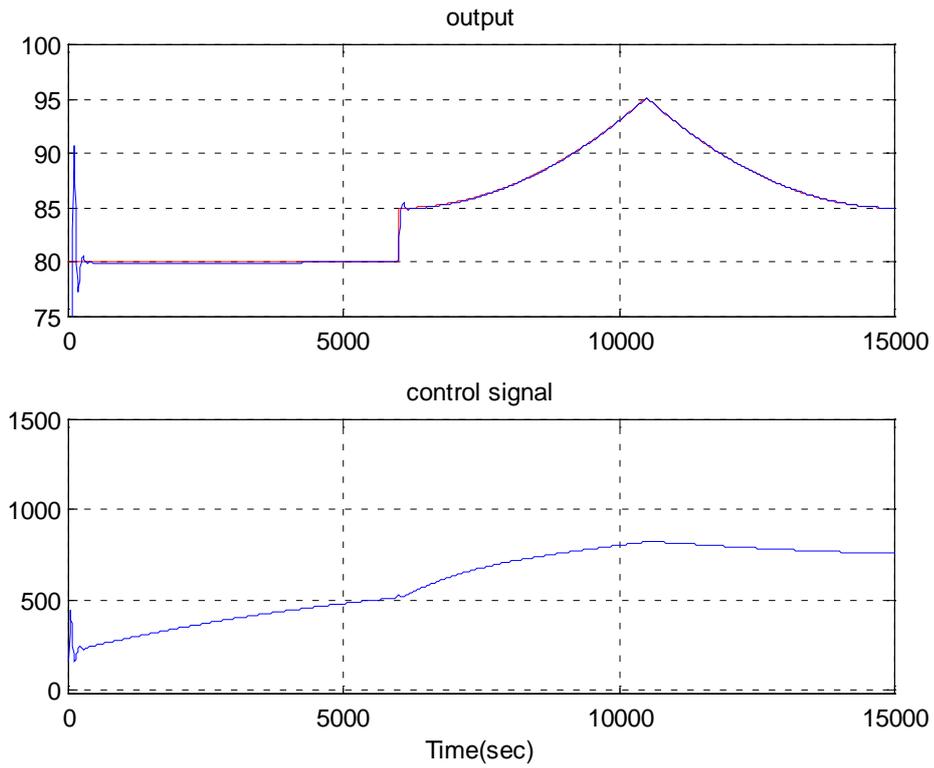

Figure 4. Controller performance in the absence of disturbance and noise

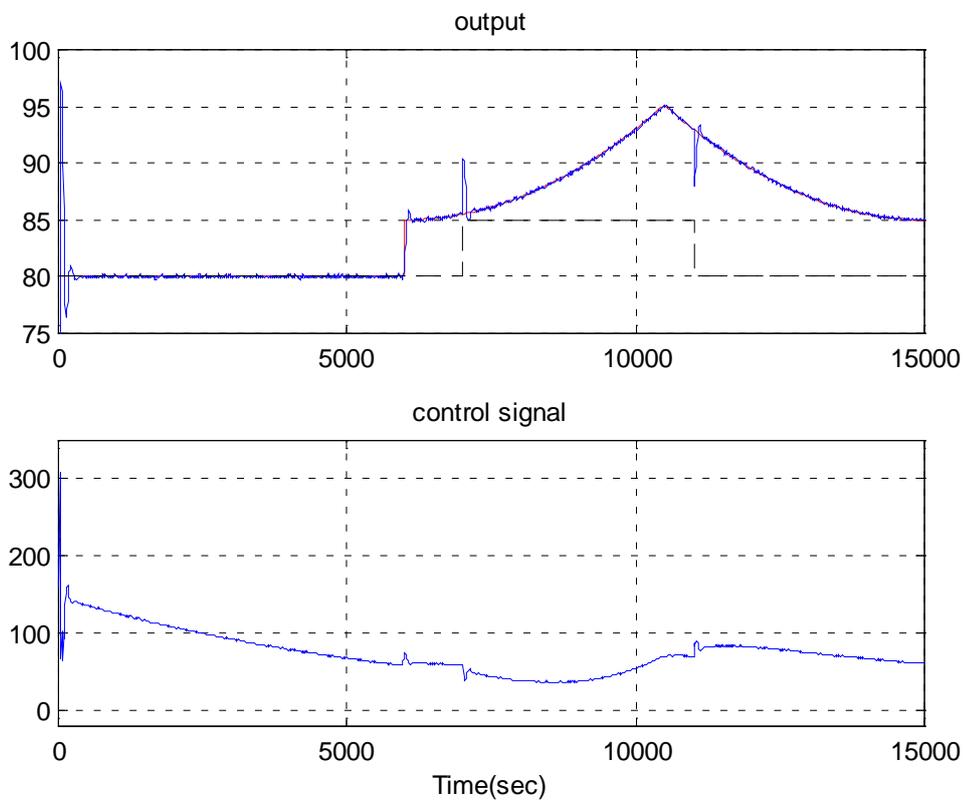

Figure 5. Controller performance in the presense of step disturbance (dashed line) and Guassian measurement noise

## 11. Conclusions

A sequential linearized model based predictive controller based on the DMC algorithm was designed to control the temperature of a batch MMA polymerization reactor. Using the mechanistic model of the polymerization, a transfer function was derived to relate the reactor temperature to the power of the heaters. The coefficients of the transfer function were calculated along the selected temperature trajectory by sequential linearization. The controller performance was studied experimentally on a small scale batch reactor.

## Acknowledgements

We would like to thank Dr. Mehdi Rafizadeh for providing the model of the reactor.